\documentclass{tac}[10pt]
\usepackage[pdftex]{graphicx}  
\usepackage[all,arc]{xy}       
\usepackage{amsmath,amssymb,latexsym}
\usepackage{amsfonts}

\newtheorem{claim}[Theorem]{Claim}

\newcommand{\bA}{{\mathbb A}}
\newcommand{\bB}{{\mathbb B}}
\newcommand{\bC}{{\mathbb C}}
\newcommand{\bD}{{\mathbb D}}
\newcommand{\bE}{{\mathbb E}}

\newcommand{\bH}{{\mathbb H}}

\newcommand{\bS}{{\mathbb S}}

\newcommand{\bY}{{\mathbb Y}}
\newcommand{\bT}{{\mathbb T}}

\newcommand{\cC}{{\mathcal C}}
\newcommand{\cD}{{\mathcal D}}
\newcommand{\cE}{{\mathcal E}}
\newcommand{\cF}{{\mathcal F}}

\newcommand{\cM}{{\mathcal M}}

\newcommand\Set{\mathbf{Set}}
\newcommand\Gpd{\mathbf{Gpd}}
\newcommand\Cat{\mathbf{Cat}}
\newcommand\Dist{\mathbf{Dist}}

\newcommand\Eq{\text{Eq}}
\newcommand\Rel{\mathbf{Rel}}
\newcommand\Span{\mathbf{Span}}

\title[Distributors and comprehensive factorization ]{Distributors and the comprehensive factorization \\ system for internal groupoids}
\author{Giuseppe Metere}
\address{via Archirafi 34, Palermo, Italy}
\eaddress{giuseppe.metere@unipa.it}
\keywords{distributor, profunctor, factorization system, internal groupoid}
\amsclass{18A32,20L05}
\thanks{This work was partially supported by the Fonds de la Recherche Scientifique -  F.N.R.S.: 2015/V 6/5/005 - IB/JN - 16385, and by the G.N.S.A.G.A.\ - Gruppo Nazionale per le Strutture Algebriche, Geometriche e le loro Applicazioni (I.N.D.A.M.).
}
\begin{document}
\maketitle
\begin{abstract}
In this note we prove that distributors between groupoids in a Barr-exact category $\cE$ form the bicategory of relations relative to the comprehensive factorization system in $\Gpd(\cE)$. The case $\cE = \Set$ is of special interest.

\end{abstract}

\section{Introduction}
Distributors (also called profunctors, or bimodules) were introduced by B\'enabou in \cite{Benabou73}.  
A fruitful approach is that of considering distributors as kind of \emph{relations between categories} (see \cite{Benabou00}, and \cite[\S7.8]{Borceux94.1}). 
In the set-theoretical case, relations can be presented as being \emph{relative} to the epi/mono factorization system.  As observed by Lawvere (\cite{Law70}), such a factorization system can be obtained from a \emph{comprehension schema}: for any set $Y$, one considers the \emph{comprehension adjunction}
\begin{equation}\label{diag:comprehension_set}
\xymatrix{
\Set/Y\ar@<1ex>[r]^-{} \ar@{}[r]|-{\bot}
&\mathbf{2}^Y\ar@<1ex>[l]
},
\end{equation}
where the category $\mathbf{2}^{Y}$ is the partially ordered set of the subsets of $Y$. For any function $\xymatrix{X\ar[r]^-f&Y}$, the (epic) unit of the adjunction provides the factorization $f=m\cdot \eta_f$, where $m$ is a mono:
$$
\xymatrix{
X\ar[dr]_-f\ar[r]^-{\eta_f}
&Im(f)\ar[d]^-m
\\&Y
}
$$

As observed in \cite{SW73}, similar arguments can be used starting with the adjunction
$$
\xymatrix{
\Cat/\bY\ar@<1ex>[r]^-{} \ar@{}[r]|-{\bot}
&\Set^\bY\ar@<1ex>[l]
},$$
but climbing one dimension up produces two distinct factorizations of a given functor:
(initial/discrete opfibration) and (final/discrete fibration). The first was named \emph{comprehensive factorization of a functor} in \cite{SW73}, as arising from a \emph{categorical comprehension schema}.

A crucial point is that the two factorization systems coincide if we consider functors between groupoids. In this note we will show that, when restricted to the category of groupoids, distributors form a bicategory of relations relative to the comprehensive factorization system. More precisely, we will prove this result in the case of internal groupoids in a Barr-exact category $\cE$; the category of groupoids is recovered for $\cE=\Set$. The key fact is the elementary observation that, for groupoids, two-sided discrete fibrations are more simply described as usual discrete fibrations (Proposition \ref{prop:two-sided}). Then, since distributors can be formulated in terms of two-sided discrete fibrations, we can relate them to the comprehensive factorization system.

The internal case is of interest for some directions of research in categorical algebra and internal category theory. For instance, concerning internal non-abelian cohomology, Bourn has developed an intrinsic version of Schreier-Mac Lane Theorem of classification of extensions using internal distributors in \cite{Bourn}, and the pointed version of a class of internal distributors, so-called \emph{butterflies}, have been studied in the semi-abelian context by Abbad, Mantovani, Metere and Vitale in \cite{AMMV} and by Cigoli and Metere in \cite{CM}. On the other hand, the non-pointed version of butterflies, called \emph{fractors} in \cite{MMV}, describe a notion of weak map between internal groupoids, where in the case of groupoids internal in groups, one recovers the notion of monoidal functor  (see \cite{Vitale}).

Finally, a description of distributor composition in terms of the associated spans was missing. With this note we aim to fill this gap, and provide a useful tool for further investigations in the area.

\section{Relations relative to a factorization system}\label{sec:Relations_in_fact}
Classically, a \emph{relation} from a set $A$ to a set $B$ is a subset $S$ of the cartesian product $A\times B$. 
For any two sets $A$ and $B$, there is a (regular) epimorphic reflection between the preorder of relations from $A$ to $B$ and the category of spans from $A$ to $B$:
$$
\xymatrix@C=8ex{
\Rel(A,B)\ \ar@<-1ex>[r]_-{i_{A,B}}\ar@{}[r]|-{\bot} &\ \Span(A,B)\ar@<-1ex>[l]_-{r_{A,B}}\,.
}
$$
The reflection is given by the (epi/mono) factorization: for a span 
$$
\xymatrix{A&\ar[l]_{e_1}E\ar[r]^{e_2}&B}
$$
one obtains its associated relation by taking the image $r_{A,B}(E)$ of the function
$$
\xymatrix@C=8ex{E\ar[r]^-{\langle e_1,e_2\rangle}&A\times B}\,.
$$
The (epi/mono) factorization system establishes also a connection between the composition of relations and the composition of spans:
  given two relations, their usual composition is precisely the reflection of their composition as spans. Globally, this means that there is  a lax biadjunction between the 2-category of relations and the bicategory of spans
$$
\xymatrix@C=10ex{
\Rel\ \ar@<-1ex>[r]_-{i}\ar@{}[r]|-{\bot} &\ \Span\ar@<-1ex>[l]_-{r}\,,
}
$$
constant on objects, where only the 2-functor $i$ is truly lax, since  $r$ is in fact a pseudo 2-functor.

More generally, one can start with any finitely complete category $\cC$ endowed with an $(\cE/\cM)$ factorization system (see Section 5.5 in \cite{Borceux94.1}). Given two objects $A$ and $B$, one defines the categories of $\cM$-relations $\Rel(A,B)$ together with the local reflections $r_{A,B} \dashv i_{A,B}$. Hence, it is possible to define the composition of $\cM$-relations as the reflection of their composition as spans, but such a composition need not be associative. As a consequence, we do not obtain a bicategory, in general. When we do get a bicategory $\Rel(\cC)$, then we call it
\begin{itemize}
\item[] the bicategory of relations in $\cC$  \emph{relative to the factorization system} $(\cE/\cM)$.
\end{itemize}
This happens, for instance, when $\cC$ is regular, or more generally, when $(\cE/\cM)$ is a proper factorization system with the class $\cE$ stable under pullbacks, but these conditions are not strictly necessary, as this article witnesses too. We will not provide further details on this general issue, but the literature on the subject is wide. The interested reader can consult \cite{Milius00} and the references therein.

\section{Internal distributors and the comprehensive factorization}
Distributors between internal categories have been introduced by B\'enabou already in \cite{Benabou73}. However, the cited reference is not as widely available as it would deserve, therefore we provide a secondary source \cite{PtJ}. 

\subsection*{Basic facts}

For the notions of internal category and internal functor in a finitely complete category $\cE$, the reader can consult \cite[B2.3]{PtJ}.  
An internal functor $F$ between internal categories $\bC$ and $\bD$ is represented by a diagram:
$$
\xymatrix@C=10ex{
C_1\ar@<-1ex>[d]_d\ar@<+1ex>[d]^c\ar[r]^-{F_1}
&D_1\ar@<-1ex>[d]_d\ar@<+1ex>[d]^c
\\
C_0\ar[u]|e\ar[r]_-{F_0}
&D_0\ar[u]|e}
$$
The functor $F$  is a \emph{discrete fibration} if and only if $c\cdot F_1=F_0\cdot c$ is a pullback. It is a \emph{discrete opfibration} if and only if $d\cdot F_1=F_0\cdot d$ is a pullback. Functors that are left orthogonal to the class of discrete fibrations are called final, those that are left orthogonal to the class of discrete opfibrations are called initial. With suitable condition on  $\cE$ (e.g.\ when $\cE$ has pullback stable reflexive coequalizers,  see \cite[Lemma B2.5.9]{PtJ}), final functors and discrete fibrations give a factorization system for the category of internal categories in $\cE$, and similarly so do initial functors and discrete opfibrations.  The latter is called \emph{comprehensive factorization}.

If $\bC$ is an internal groupoid in $\cE$, the internal inverse map is denoted by $\tau\colon C_1\to C_1$.  
Bourn has shown in \cite {Bourn} that, if the base category $\cE$ is Barr-exact, then the category $\Gpd(\cE)$ of groupoids in $\cE$ admits the comprehensive factorization. Notice that, in this case, discrete fibrations coincide with discrete opfibrations, and final functors with initial functors. Therefore, the two factorization systems mentioned above reduce to a single one that we denote by $(\cF/\cD)$. From now on, we will assume that $\cE$ is Barr-exact.

The \emph{connected components} functor
$$
\xymatrix{\Pi_0\colon\Gpd(\cE)\ar[r]& \cE}
$$
is defined: it assigns to every internal groupoid, the coequalizer in $\cE$ of its domain and codomain maps. Recall that since $\cE$ is Barr-exact, groupoids in $\cE$  have effective support, i.e.\ the regular image of the map $\langle d,c \rangle$ coincides with the kernel pair relation of such coequalizer. 
Cigoli in \cite{Cigoli2017} has characterized final functors between groupoids in a Barr-exact category $\cE$. 
\proposition[\cite{Cigoli2017}]\label{prop:Alan}
An internal functor $F\colon \bC\to\bD$ between groupoids in a Barr-exact category $\cE$ is  final if and only if it is internally full and essentially surjective, i.e.\ if and only if 
\begin{itemize}
\item the canonical comparison of $C_1$ with the joint pullback of $d$ and $c$ along $F_0$ is a regular epimorphism;
\item $\Pi_0(F)$ is a regular epimorphism.
\end{itemize}
\endproposition
Let us notice that, if $F$ is a full functor, then the morphism $\Pi_0(F)$ is a monomorphism, so that in the proposition above, it is an isomorphism.

\subsection*{Distributors between groupoids are discrete fibrations}
The definition of internal distributor closely follows the set-theoretical definition. 
\definition\label{def:relator} (\cite{PtJ})
Let $\bA$ and $\bB$ be internal groupoids in $\cE$. A distributor
$$
\xymatrix{\bB&\ar[l]|-@{|}_-{S}\bA}  
$$
consists of the following data:
\begin{itemize}
\item a span $\xymatrix{A_0&\ar[l]_{L}S_0\ar[r]^{R}&B_0}$ in $\cE$,
\item a left action $\xymatrix{A_1\underset{A_0}{\times} S_0\ar[r]^-{\lambda_S}&S_0}$\,,
\item a right action $\xymatrix{S_0\underset{B_0}{\times}B_1\ar[r]^-{\rho_S}&S_0}$\,,
\end{itemize}
which are associative, unital and compatible, where compatible means that the following diagram commutes:
$$
\xymatrix@C=10ex{
A_1\underset{A_0}{\times}S_0\underset{B_0}{\times}B_1\ar[r]^-{1\times\rho_S}\ar[d]_{\lambda_S\times 1} 
&A_1\underset{A_0}{\times}S_0\ar[d]^{\lambda_S}
\\
S_0\underset{B_0}{\times}B_1\ar[r]_-{\rho_S}
&S_0
}
$$ 
\enddefinition
Distributors between two given groupoids $\bA$ and $\bB$ form the category $\Dist(\bA,\bB)$, where an arrow between two distributors
$$
\alpha\colon (L,S_0,R)\to (L',S_0',R')
$$
is an arrow in the base category $\alpha\colon S_0\to S_0'$ such that $L'\cdot\alpha=L$,  $R'\cdot\alpha=R$, commuting with the actions.

Like in the set-theoretical case, every internal distributor determines a span in $\Gpd(\cE)$. For instance, the distributor $S$ above determines the span
\begin{equation}\label{diag:span}
\xymatrix{\bA&\ar[l]_{L}\bS\ar[r]^{R}&\bC}
\end{equation}
where the internal groupoid $\bS$ has $S_0$ as the object of objects, and the object of arrows $S_1$  is obtained by the pullback
$$
\xymatrix@C=10ex{
S_1\ar[r]^{\pi_2}\ar[d]_{\pi_1}
&S_0\underset{B_0}{\times}B_1\ar[d]^{\rho_S}
\\
A_1\underset{A_0}{\times}S_0\ar[r]_-{\lambda_S}
&S_0
}
$$
with structure maps 
$$
\xymatrix@!=8ex{
d\colon S_1\ar[r]^-{\pi_1}
&A_1\underset{A_0}{\times}S_0\ar[r]^-{\pi_2}
&S_0
}
$$
$$
\xymatrix@!=8ex{
c\colon S_1\ar[r]^-{\pi_2}
&S_0\underset{B_0}{\times}B_1\ar[r]^-{\pi_1}
&S_0
}
$$ 
and $e\colon S_0\to S_1$ is the unique morphism such that $\pi_1\cdot e=\langle e\cdot L_0,1 \rangle$ and $\pi_2\cdot e = \langle 1,e\cdot R_0 \rangle$.
Finally, the internal functors $L$ and $R$ are described below:
$$
\xymatrix{
S_1\ar@<-1ex>[d]_d\ar@<+1ex>[d]^c\ar[r]^-{\pi_1}
&A_1\underset{A_0}{\times}S_0\ar[r]^-{\pi_1}
&A_1\ar@<-1ex>[d]_d\ar@<+1ex>[d]^c
\\S_0\ar[u]|e\ar[rr]_{L_0}&&A_0\ar[u]|e
}
\qquad
\xymatrix{
S_1\ar@<-1ex>[d]_d\ar@<+1ex>[d]^c\ar[r]^-{\pi_2}
&S_0\underset{B_0}{\times}B_1\ar[r]^-{\pi_2}
&A_1\ar@<-1ex>[d]_d\ar@<+1ex>[d]^c
\\S_0\ar[u]|e\ar[rr]_{R_0}&&B_0\ar[u]|e
}
$$

The following result establishes the connection between the notion of distributor between  groupoids and that of discrete fibration. 
\proposition\label{prop:two-sided}
Giving a distributor
$$
\xymatrix{\bB&\ar[l]|-@{|}_-{S}\bA}  
$$
between internal groupoids in a finitely complete category is equivalent to giving an internal discrete fibration 
\begin{equation}\label{diag:prd}
\xymatrix@C=10ex{\bS\ar[r]^-{\langle L,R\rangle}&\bA\times\bB}
\end{equation}
\endproposition
\proof
A span  $(L,R)$ as in  (\ref{diag:span}) is  determined by a distributor if and only if it is a two-sided discrete fibration.  For groupoids in $\Set$, it is easy to prove that such a span  is a two-sided discrete fibration if and only if the induced functor into the product (\ref{diag:prd})
is a discrete fibration. The result for internal groupoids follows by the usual Yoneda embedding argument.
\endproof

\remark
The notion of two-sided discrete fibration, introduced by Street in \cite{Street}, appears earlier in the literature, although implicitly, as a discretization of so-called regular spans, introduced and studied by  Yoneda in \cite{Yo60} (see also \cite{CMMV18}). It is relevant to our discussion to recall that Yoneda, in Section 3.5 of the cited paper, introduces a (generalized) composition product of two composable regular spans as a suitable discretization of their composition as spans.  This was the starting point for our investigations on the subject. 
\endremark

The last proposition allows us to describe the reflection of spans into distributors. Since $(\cF/\cD)$ is a factorization system, we need not prove the following statement.
\proposition\label{prop:reflect_span_in_dist}
Let $\bA$ and $\bB$ be two groupoids in a Barr-exact category $\cE$. The comprehensive factorization defines the reflection $R$ of spans between $\bA$ and $\bB$ into distributors between $\bA$ and $\bB$, i.e.\ it establishes an adjoint pair
$$
\xymatrix{
\Dist(\bA,\bB) \ar@<-1ex>[r]_{I_{\bA,\bB}} \ar@{}[r]|{\bot}
&\Span(\bA,\bB) \ar@<-1ex>[l]_{R_{\bA,\bB}}
}
$$
with $R_{\bA,\bB}\cdot I_{\bA,\bB} \simeq id_{\Dist(\bA,\bB)}$.
For the span $(P,Q)$, the reflection is obtained  by the factorization
$$
\xymatrix{
\bE\ar[dr]_{\langle P,Q\rangle}\ar[r]^-{F}
&R(\bE)\ar[d]^{\langle L,R\rangle}
\\&\bA\times\bB
}
$$
where $F$ is final, and  $\langle L,R\rangle$ a discrete fibration.
\endproposition

\subsection*{Composition of distributors}

In this section we assume the base category $\cE$ to be Barr-exact.
With this hypothesis, internal groupoids in $\cE$ admit $\Pi_0$'s, and these are stable under pullback. Therefore, internal distributors in $\cE$ can be composed.

For two composable distributors $T$ and $S$
$$
\xymatrix{\bC&\ar[l]|-@{|}_-{T}\bB&\ar[l]|-@{|}_-{S}\bA}  
$$
their composition $\xymatrix{A_0&\ar[l]_-{L}(T\otimes S)_0\ar[r]^-{R}&C_0}$ is obtained by the universal property of the coequalizer $Q_0$ in the first line in the diagram below
\begin{equation}\label{diag:coequalizer}
\begin{aligned}
\xymatrix@C=10ex{
S_0\underset{B_0}{\times}B_1\underset{B_0}{\times}T_0\ar@<+1ex>[r]^-{\rho_S\times 1} \ar@<-1ex>[r]_-{1\times\lambda_T} 
&S_0\underset{B_0}{\times}T_0\ar[r]^-{Q_0}\ar[dr]\ar@{.>}[l]
&(T\otimes S)_0\ar@{-->}[d]^{\exists !}
\\&&A_0\times C_0
}
\end{aligned}
\end{equation}
where, in order to simplify notation, we write $S_0\underset{B_0}{\times}B_1\underset{B_0}{\times}T_0$ for the limit $S_0\,_{R} {\times}_{c}\,B_1 \,_{d} {\times}_L\,T_0$ and $S_0\underset{B_0}{\times}T_0$ for the pullback $S_0\,_{R}{\times}_L\,T_0$.
Thanks to pullback stability of $\Pi_0's$, the actions $\rho_{T\otimes S}$ and $\lambda_{T\otimes S}$ are induced by $\rho_T$ and $\lambda_S$, see \cite[B2.7]{PtJ} for details.
\remark\label{rk}
The first line of diagram (\ref{diag:coequalizer}) can be interpreted as a groupoid $\bH$, together with its object of connected components $\Pi_0(\bH)=(T\otimes S)_0$. For the reader's convenience,  we provide a set-theoretical description of such a groupoid. The objects of $\bH$,  i.e.\ the elements of $S_0\underset{B_0}{\times}T_0$, are pairs $(s,t)$ such that $R(s)=L(t)$. The elements of  $S_0\underset{B_0}{\times}B_1\underset{B_0}{\times}T_0$ are the arrows of $\bH$. More precisely, one such element $(s,\beta,t)$ is an arrrow 
$(s\cdot \beta, t)\to (s,\beta\cdot t)$, where we have used the right and the left action of $\bB$ on $S_0$ and $T_0$ respectively. Arrows composition is inherited from arrows composition in $\bB$. Under this interpretation, it is clear that the coequalizer $(T\otimes S)_0$ represents the connected components of $\bH$.
\endremark
The next statement is the key result of this note. It relates distributor composition to span composition.
\proposition\label{prop:composition_agrees}
Distributor composition agrees with (the reflection of) span composition, i.e\  for given groupoids $\bA$, $\bB$ and $\bC$, the following diagram commutes:
$$
\xymatrix{
\Dist(\bA,\bB)\times\Dist(\bB,\bC)\ar[r]^-{\otimes}\ar[d]_-{I\times I}
&\Dist(\bA,\bC)
\\
\Span(\bA,\bB)\times\Span(\bB,\bC)\ar[r]_-{\diamond}
&\Span(\bA,\bC)\ar[u]_-{R}
}
$$
where $\otimes$ is the composition of distributors and $\diamond$ is the composition of spans. 
\endproposition
\proof
The way the composite $R\cdot \diamond \cdot I\!\times\!I$ acts on a pair of distributors $S$ and $T$ is shown in the following diagram
$$
\xymatrix@!=3ex{
&&T\diamond S\ar[dr]|!{[ddll];[rr]}\hole^{\bar R} \ar[dl]_{\bar L} \ar[rr]^-{F}
&&R(T\diamond S)\ar[ddllll]_{\hat L}\ar[dd]^{\hat R}
\\
&\bS\ar[dr]|!{[rrru];[dl]}\hole_(.6){R} \ar[dl]_{L}&&\bT\ar[dr]_{R}\ar[dl]^{L}\\
\bA&&\bB&&\bC
}
$$
where the square $R\cdot\bar L =L\cdot \bar R$ is a pullback in $\Gpd(\cE)$.
Hence, we consider two factorizations of the functor $\langle L\cdot \bar L, R\cdot \bar R\rangle$:
\begin{equation}\label{diag:Q}
\begin{aligned}
\xymatrix@C=20ex{
T\diamond S\ar[r]^-{F}\ar[dr]|{\langle L\cdot \bar L, R\cdot \bar R\rangle} \ar[d]_{Q}
&R(T\diamond S)\ar[d]^{\langle \hat L,\hat R\rangle}
\\
T\otimes S\ar[r]_-{\langle L,R\rangle}
&\bA\times \bC
}
\end{aligned}
\end{equation}
The first one is the comprehensive factorization, that is given by the final functor $F$ followed by the discrete fibration $\langle \hat L,\hat R\rangle$. The second one extends at the arrows level the factorization provided by the coequalizer diagram in (\ref{diag:coequalizer}). It consists of a functor $Q$ (description below) followed by the discrete fibration  $\langle L,R\rangle$ representing the distributor $T\otimes S$. By uniqueness of factorization, it suffices to prove that $Q$ is final in order to prove that these two factorizations are isomorphic. 

First we need to recall that the pullback groupoid $T\diamond S$ is computed levelwise in $\cE$, therefore, following the lines of the simplified notation adopted in diagram (\ref{diag:coequalizer}) above, we have $(T\diamond S)_i=S_i\underset{B_i}{\times}T_i$, for $i=0,1$. Moreover domain, codomain and unit maps are given by the universal properties of such pullbacks, namely $\langle d,d\rangle$, $\langle c,c \rangle$ and $\langle e,e\rangle$ respectively. 

We are ready to describe the functor $Q$ explicitly. For internal groupoids, the internal two-sided discrete fibration associated with  $T\otimes S$ is just a discrete fibration, and the cited factorization can be represented as follows:
\begin{equation}\label{diag:internal_fact}
\begin{aligned}
\xymatrix@C=10ex{
(T\diamond S)_1 \ar@<-1ex>[d]_{\langle d,d \rangle} \ar@<+1ex>[d]^{\langle c,c \rangle} \ar@{-->}[r]^-{\exists ! Q_1}
&(T\otimes S)_0\underset{A_0\!\times C_0}{\times}(A_1\times C_1)\ar@<-1ex>[d]_{\bar d} \ar@<+1ex>[d]^{\bar c} \ar[r]^-{\pi_2}
&A_1\times C_1\ar@<-1ex>[d]_{d\times d} \ar@<+1ex>[d]^{c\times c}
\\
(T\diamond S)_0\ar[r]_-{Q_0}\ar[u]
&(T\otimes S)_0\ar[r]_-{\langle L_0,R_0\rangle}\ar[u]
&A_0\times C_0\ar[u]}
\end{aligned}
\end{equation}
where the downward directed squares on the right are pullbacks.
 By Proposition \ref{prop:Alan}, $Q=(Q_1,Q_0)$ is final if and only if it is  internally full and essentially surjective. 
Therefore, the proof of the proposition will be achieved through the proof of the following two claims.
\endproof
\claim\label{claim2} The arrow 
$$
\xymatrix@C=12ex{\Pi_0(Q)\colon\Pi_0(T\diamond S)\ar[r]&\Pi_0(T\otimes S)}
$$
is a regular epimorpism.
\endclaim

\proof[of Claim \ref{claim2}]
The arrow $\Pi_0(Q)$ is a regular epimorphism since $Q_0$ is. 
\endproof
\claim\label{claim1} The comparison map $K$ with the joint pullback $W$ in the diagram below is a regular epimorphism. 
\begin{equation}\label{diag:K}
\begin{aligned}
\xymatrix{
(T\diamond S)_1\ar@{-->}[dr]^{K}\ar@/^3ex/[drr]^{Q_1}\ar@/_2ex/[ddr]_{\langle\langle c,c \rangle,\langle d,d \rangle\rangle}
\\
&W\ar[d]_-{J}\ar[r]^-{K'}
&(T\otimes S)_0\underset{A_0\!\times C_0}{\times}(A_1\times C_1)\ar[d]^{\langle \bar c,\bar d\rangle}
\\
&(T\diamond S)_0\times(T\diamond S)_0\ar[r]_{Q_0\times Q_0}
&(T\otimes S)_0\times (T\otimes S)_0
}
\end{aligned}
\end{equation}
\endclaim

\proof[of Claim \ref{claim1}]
It is convenient to describe  $W$ by means of the pasting of the two pullbacks below:
\begin{equation}\label{diag:K2}
\begin{aligned}
\xymatrix@!=7ex{
&&W\ar[dl]_{K''}\ar[ddrr]^{J_2}
\\
&(T\diamond S)_0\underset{A_0\!\times C_0}{\times}(A_1\times C_1) \ar[dl]_{\pi_1}\ar[dr]^{Q_0\times 1}
\\
(T\diamond S)_0\ar[dr]_{Q_0}
&&(T\otimes S)_0\underset{A_0\!\times C_0}{\times}(A_1\times C_1) \ar[dl]^{\bar c}\ar[dr]_{\bar d}
&&(T\diamond S)_0\ar[dl]^{Q_0}
\\
&(T\otimes S)_0&&(T\otimes S)_0
}
\end{aligned}
\end{equation}
with $J=\langle \pi_1\cdot K'',J_2\rangle$ and $K'=(Q_0\times 1) \cdot K''$. 

Now,  if we denote by $\lambda_{T\otimes S}$ the left action and by $\rho_{T\otimes S}$ the right action associated with the distributor composition of $S$ and $T$, the compatibility axiom for these two action is expressed by the the equality of the two morphisms
$\rho_{T\otimes S}\cdot(\lambda_{T\otimes S}\times 1)$ and $\lambda_{T\otimes S} \cdot (1\times \rho_{T\otimes S})$. Let us denote these equal morphisms by $\lambda_{S}\otimes\rho_{T}$, and then consider the commutative diagram:
\begin{equation}\label{diag:D}
\begin{aligned}
\xymatrix@C=14ex{
(T\diamond S)_0\underset{A_0\!\times C_0}{\times}(A_1\times C_1)\ar[r]^{Q_0\times 1}\ar[d]^{\text{iso}}\ar@/_12ex/[dd]_{D}
&(T\otimes S)_0\underset{A_0\!\times C_0}{\times}(A_1\times C_1)\ar@/^12ex/[dd]^{\bar d}\ar[d]_{\text{iso}}
\\
A_{1}\underset{A_0}{\times}S_{0}\underset{B_0}{\times}T_{0}\underset{C_0}{\times}C_{1}\ar[d]^{\lambda_{S}\times\rho_{T}}\ar[r]^{1\times Q_{0}\times 1}
&A_{1}\underset{A_0}{\times}(T\otimes S)_{0}\underset{C_0}{\times}C_{1}\ar[d]_{\lambda_{S}\otimes\rho_{T}}
\\
(T\diamond S)_0\ar[r]_{Q_0}
&(T\otimes S)_0
}
\end{aligned}
\end{equation}
where the isomorphisms above are clear (notice that $D$ is defined by the diagram).
Using the equation $\bar d\cdot (Q_0\times 1)= Q_0 \cdot D$, we can obtain a  description of $W$ through the kernel pair of $Q_0$. 
This is shown in the diagram below, which deserves some explanation.
\begin{equation}\label{diag:K3}
\begin{aligned}
\xymatrix@!=7ex{
&&(T\diamond S)_1\ar[dl]_K \ar[dr]^M\ar@{}[dd]|{(I)}
\\
&W\ar[dl]_{K''} \ar[dr]^H \ar@{}[dd]|(.27){(II)}
&&S_0\underset{B_0}{\times}B_1\underset{B_0}{\times}T_0\ar[dl]^{\Sigma}
\\
(T\diamond S)_0\underset{A_0\!\times C_0}{\times}(A_1\times C_1)\ar[dr]_D
&&\Eq(Q_0)\ar[dl]_{\pi_1}\ar[dr]^{\pi_2} \ar@{}[dd]|{(III)}
\\
&(T\diamond S)_0\ar[dr]_{Q_0}
&&(T\diamond S)_0\ar[dl]^{Q_0}
\\
&&(T\otimes S)_0
}
\end{aligned}
\end{equation}
The arrow $H$ is determined by the equations $\pi_1\cdot H=D\cdot K''$ and $\pi_2\cdot H=J_2$. The pasting of $(II)$ and $(III)$ in (\ref{diag:K3}) corresponds to the big rectangular pullback in diagram (\ref{diag:K2}). Therefore, since $(III)$ is a pullback, so is $(II)$.
In order to continue the description of the diagram, we let the arrow $\Sigma$ be the comparison with the kernel pair $(\Eq(Q_0),\pi_1,\pi_2)$ of $Q_0$, i.e.\ the unique arrow satisfying the equations $\pi_1\cdot \Sigma = \rho_S\times 1$ and $\pi_2\cdot \Sigma = 1\times \lambda_T$.
Here comes a main point where Barr-exactness of the base category (as opposed to regularity) is actually used. 
In fact, by exactness, the support of the groupoid $\bH$
(see Remark \ref{rk}) is effective, i.e.\ it coincides with the kernel pair relation of the coequalizer of  domain and codomain of  $\bH$. This means that the comparison $\Sigma$ is the regular epimorphic part of the joint  factorization of domain and codomain of  $\bH$ through its support.

Finally, the construction of the arrow $M$ requires some more work. Let us consider the following diagram:
$$
\xymatrix@C=12ex{
S_0\underset{B_0}{\times}B_1\underset{B_0}{\times}T_0 \ar[d]_{\rho_S\times 1} \ar[r]^-{\pi_2}   
&B_1\ar[d]^d
\\
(T\diamond S)_0\ar[r]_{R_0\cdot\bar L_0=L_0\cdot\bar R_0}
&B_0
}
$$
one easily checks that it is a pullback. Then, the arrow $M$ is uniquely determined by the equations 
$$
(\rho_S\times 1)\cdot M = D\cdot K''\cdot K\,,\qquad\pi_2\cdot M = \pi_{B_1}
$$
where $\pi_{B_1}$ is the diagonal of the pullback defining $(T\diamond S)_1$.
At this point it is possible to prove  that $(I)$ commutes by composing with pullback projections. First one easily compute $\pi_1\cdot\Sigma\cdot M= (\rho_S\times 1)\cdot M =D\cdot K''\cdot K=\pi_1\cdot H\cdot K$.
Then, with a little more effort,  $\pi_2\cdot\Sigma\cdot M= (1\times \lambda _T)\cdot M =\langle c,c\rangle= J_2\cdot K=\pi_2 \cdot H \cdot K$, where the only non-obvious equality is the second one. However, it can be easily proved with elements, and then use Yoneda embedding.

Yoneda embedding can also be used to prove that $(I)+(II)$ is a pullback. Therefore, $(I)$ is a pullback too, and since $\Sigma$ is a regular epimorphism, so is $K$.
\endproof

\remark
By Corollary 2.5 in \cite{Cigoli2017},  since $Q$ is full, $\Pi_0(Q)$ is also a monomorphism, and therefore it is an isomorphism.
\endremark

We conclude the section with the expected result.

\theorem\label{thm:main2}
Let $\cE$ be a Barr-exact category. Then
$$
\Dist{\Gpd}(\cE)=\Rel(\Gpd(\cE))\ {\text w.r.t.}\  (\cF/\cD),
$$
i.e.\ the bicategory of distributors between internal groupoids in $\cE$ is the bicategory of relations in $\Gpd(\cE)$ relative to the (final/discrete fibration) factorization system.
\endtheorem
\begin{proof}
Proposition \ref{prop:reflect_span_in_dist} identifies $\cD$-relations, and Proposition \ref{prop:composition_agrees} provides the bicategory structure.
\end{proof}

\section*{Acknowledgments}
I am grateful to Alan Cigoli and to Sandra Mantovani for the many interesting conversations and suggestions, and to Tim Van der Linden for supporting my visit in Louvain-la-Neuve in 2016. 

\smallskip
The first version of this paper contains mistakes, I am in debt to Steve Lack for pointing this out to me.  A critical fact is that Lemma 3.2 of first version is not correct. In the present version, we provide a different proof of Proposition \ref{prop:composition_agrees} which does not use the cited lemma, and we fix the other issues.

\end{document}